\documentclass[a4paper,pdftex,reqno]{amsart}

\makeatletter
\@namedef{subjclassname@2020}{\textup{2020} Mathematics Subject Classification}
\makeatother

\usepackage{geometry}
\title[Twisted Alexander invariants for the braid group]%
{Twisted Alexander invariants for the braid group associated with the Tong-Yang-Ma representation}
\author{Akihiro Takano}


\address{GRADUATE SCHOOL OF MATHEMATICAL SCIENCES, THE UNIVERSITY OF TOKYO, 3-8-1 KOMABA, MEGURO-KU, TOKYO, 153-8914, JAPAN}
\email{takano@ms.u-tokyo.ac.jp}

\usepackage{amsmath} 
\usepackage{amsfonts}
\usepackage{mathrsfs}
\usepackage{amssymb}
\usepackage{amsthm}
\usepackage{bm}
\usepackage[new]{old-arrows}
\usepackage{wrapfig}
\usepackage{graphicx}
\usepackage{color}
\usepackage[labelsep=colon]{caption}
\usepackage[labelformat=empty,subrefformat=parens]{subcaption}
\usepackage{multicol}
\usepackage{units}
\usepackage{hyperref}

\newcommand{\relmiddle}[1]{\mathrel{}\middle#1\mathrel{}}

\newtheorem{thm}{Theorem}[section]

\newtheorem{lemma}[thm]{Lemma}
\newtheorem{cor}[thm]{Corollaly}

\theoremstyle{definition}
\newtheorem{defi}[thm]{Definition}

\newtheorem{q}[thm]{Question}

\newcommand{\Z}{\mathbb{Z}}

\newcommand{\wtil}{\widetilde}
\newcommand{\dfr}{\displaystyle\frac}
\newcommand{\B}{\mathscr{B}}

\begin{document}

\begin{abstract}
In this paper, we compute the twisted Alexander invariant of the braid group associated with the Tong-Yang-Ma representation.
\end{abstract}

\maketitle

\section{Introduction}
The twisted Alexander invariant for a finitely presentable group was introduced by Wada \cite{Wada}.
This is an invariant for the given group and its representation.
The twisted Alexander invariant is a generalization of the Alexander invariant and in paritular, if we take the trivial representation, then the twisted Alexander invariant almost coincides with the Alexander invariant.

Morifuji \cite{Morifuji} studied the twisted Alexander invariants associated with the Jones representations.
In particular, he calculated the invariant associated with the reduced Burau representation and its result is the following:

\begin{thm} [{\cite[Theorem 1.1]{Morifuji}}]
Let $\wtil{\B}_n \colon B_n \longrightarrow GL_{n-1}(\mathbb{Z}[t^{\pm1}])$ be the reduced Burau representation of the braid group $B_n$ and $\alpha \colon B_n \longrightarrow \mathbb{Z} \cong \langle z \rangle$ the abelianization. Then the twisted Alexander invariant $\Delta_{B_n, \wtil{\B}_n}(z)$ is given by
$$
\Delta_{B_n, \wtil{\B}_n}(z)=
\begin{cases}
1-tz^2 & (n=3)\\
1 & (n \geq 4)
\end{cases}.
$$
\end{thm}

Tong, Yang, and Ma \cite{Tong-Yang-Ma} researched the representations of the braid group $B_n$ such that the $i$-th generator in the Artin presentation maps to the regular matrix $I_{i-1} \oplus T \oplus I_{n-i-1}$ where $I_{k}$ is the $k \times k$ identity matrix and $T$ is an $m \times m$ regular matrix whose entries are elements of $\mathbb{Z} [t^{\pm1}]$.
They proved that there exist three kinds of irreducible representations:
the trivial one, the Burau one, and a new $n$-dimensional one.
In the case of $m=2$, there essentially exist only two non-trivial representations:
one is the unreduced Burau representation, and the other is the irreducible representation called the Tong-Yang-Ma representation.

In this paper, we show that the twisted Alexander invariant of the braid group associated with the Tong-Yang-Ma representation is similar to the above.
More precisely, we have the following:

\begin{thm} \label{twisted_TYM}
Let $TYM_n \colon B_n \longrightarrow GL_{n}(\mathbb{Z}[t^{\pm1}])$ be the Tong-Yang-Ma representation of the braid group $B_n$ and $\alpha \colon B_n \longrightarrow \mathbb{Z} \cong \langle z \rangle$ the abelianization. Then the twisted Alexander invariant $\Delta_{B_n,TYM_n}(z)$ is given by
$$
\Delta_{B_n,TYM_n}(z)=
\begin{cases}
1+tz^3 & (n=3)\\
1 & (n \geq 4)
\end{cases}.
$$
\end{thm}

Moreover, there is an extension of the braid group called the welded braid group $WB_n$, and there are some representations of $WB_n$ extended from representations of the braid group (see \cite{Bellingeri-Soulie}).
The Tong-Yang-Ma representation is one of these representations.
Thus, we also compute the twisted Alexander invariant of the welded braid group associated with the Tong-Yang-Ma representation for $n=3$.

\section{Tong-Yang-Ma representation} \label{tym}
Let $B_n$ be the braid group of $n$ strings. The group $B_n$ has the following presentation, which is called the Artin presentation:
$$
\left\langle \sigma_1 , \ldots , \sigma_{n-1} \relmiddle|
\begin{array}{ll}
\sigma_i \sigma_j = \sigma_j \sigma_i & (|i-j| \geq 2) \\ 
\sigma_i \sigma_{i+1} \sigma_i = \sigma_{i+1} \sigma_i \sigma_{i+1} & (i=1, \ldots ,n-2)
\end{array}
\right\rangle.
$$
Tong, Yang, and Ma \cite{Tong-Yang-Ma} studied how many representations $B_n \longrightarrow GL_{n}(\mathbb{Z} [t^{\pm1}])$ there are of the form
$$
\sigma_i \longmapsto I_{i-1} \oplus
\left(
\begin{array}{cc}
a & b \\
c & d
\end{array}
\right)
\oplus I_{n-i-1}
$$
for $i = 1, \ldots, n-1$.
They found that there essentially, namely up to equivalent, transposition and constant multiplied, exist only two non-trivial representations of this type.
One is the unreduced Burau representation $\B_n \colon B_n \longrightarrow GL_n(\mathbb{Z}[t^{\pm1}])$, that is
$$
\sigma_i \longmapsto
I_{i-1} \oplus
\left(
\begin{array}{cc}
0 & t \\
1 & 1-t
\end{array}
\right)
\oplus I_{n-i-1},
$$
and the other is the irreducible representation given by
$$
\sigma_i \longmapsto
I_{i-1} \oplus
\left(
\begin{array}{cc}
0 & 1 \\
t & 0
\end{array}
\right)
\oplus I_{n-i-1}.
$$
The latter representation is called the \textbf{Tong-Yang-Ma representation}
$$
TYM_n \colon B_n \longrightarrow GL_{n} (\mathbb{Z} [t^{\pm1}]).
$$

\section{Twisted Alexander invariant}
In this section, we refer to Wada \cite{Wada}.
Let $G$ be a group with a finite presentation
\begin{align} \label{pre}
G=\langle x_1 , \ldots , x_l \mid r_1 , \ldots , r_m \rangle .
\end{align}
Suppose that $G$ has a surjective homomorphism $\alpha \colon G \longrightarrow \mathbb{Z} \cong \langle z \rangle$.
Let $\rho \colon G \longrightarrow GL_n(R)$ be a linear representation, where $R$ is a unique factorization domain.
Extending these maps linearly to the group ring $\mathbb{Z}[G]$, we obtain two ring homomorphisms
$$
\widetilde{\alpha} \colon \mathbb{Z}[G] \longrightarrow \mathbb{Z}[z^{\pm1}]\ \ \textrm{and}\ \ \widetilde{\rho} \colon \mathbb{Z}[G] \longrightarrow M_n(R),
$$
where $M_n(R)$ is the matrix algebra of degree $n$ over $R$.
Then the tensor product homomorphism $\widetilde{\rho} \otimes \widetilde{\alpha} \colon \mathbb{Z}[G] \longrightarrow M_n(R[z^{\pm1}])$ of $\widetilde{\rho}$ and $\widetilde{\alpha}$ is defined by
$$
(\widetilde{\rho} \otimes \widetilde{\alpha}) (g) := \widetilde{\rho} (g) \widetilde{\alpha} (g)\ \ \ (g \in G).
$$
Let $F_l$ be the free group generated by $x_1, \ldots, x_l$ and $\phi \colon F_l \longrightarrow G$ the surjective homomorphism induced by each presentation.
Similarly to $\alpha$ and $\rho$, $\phi$ induces a ring homomorphism
$$
\widetilde{\phi} \colon \mathbb{Z}[F_l] \longrightarrow \mathbb{Z}[G].
$$
Then the composition map $\Phi := (\widetilde{\rho} \otimes \widetilde{\alpha}) \circ \widetilde{\phi} \colon \mathbb{Z}[F_l] \longrightarrow M_n(R[z^{\pm1}])$ is a ring homomorphism.

We define the $m \times l$ matrix $M$ whose $(i,j)$ component is the $n \times n$ matrix
$$
\Phi \left( \frac{\partial r_i}{\partial x_j} \right) \in M_n(R[z^{\pm1}]),
$$
where $\displaystyle{\frac{\partial}{\partial x_j}}$ is the Fox derivation with respect to $x_j$, that is, a $\Z$-linear map $\mathbb{Z}[F_l] \longrightarrow \mathbb{Z}[F_l]$ satisfying two conditions:
\begin{itemize}
\item $\dfr{\partial x_i}{\partial x_j} = \delta_{ij}$, where $\delta_{ij}$ is the Kronecker delta, and
\item $\dfr{\partial (g g')}{\partial x_j} = \dfr{\partial g}{\partial x_j} + g \dfr{\partial g'}{\partial x_j}$ for any $g, g' \in F_l$.
\end{itemize}
This matrix $M$ is called the \textbf{Alexander matrix} of the presentation (\ref{pre}) of $G$ associated with the representation $\rho$.

For $1 \leq j \leq l$, let $M_j$ be the $m \times (l-1)$ matrix obtained from $M$ by removing the $j$-th column. We regard $M_j$ as an $mn \times (l-1)n$ matrix with coefficients in $R[z^{\pm1}]$.
For an $(l-1)n$-tuple of indices
$$
I=\left( i_1 , \ldots , i_{(l-1)n} \right) \ \ \ \left(1 \leq i_1 < \cdots < i_{(l-1)n} \leq mn \right),
$$
we write $M^I_j$ for the $(l-1)n \times (l-1)n$ square matrix consisting of the $i_k$-th rows of the matrix $M_j$, where $k=1 , \ldots , (l-1)n$.

In order to define the twisted Alexander invariant, we prepare the following two lemmas:

\begin{lemma} [{\cite[Lemma 2]{Wada}}] \label{twisted1}
There is an integer $j\ (1 \leq j \leq l)$ such that $\det \Phi (1-x_j) \neq 0$.
\end{lemma}

\begin{lemma} [{\cite[Lemma 3]{Wada}}] \label{twisted2}
For any integers $j, k\ (1 \leq j < k \leq l)$ and choice of indices $I$,
$$
\det M^I_j \det \Phi (1-x_k) = \pm \det M^I_k \det \Phi (1-x_j).
$$
In fact, if the dimension of the representation $\rho$ is even, then the sign in this formula is $+$.
\end{lemma}

By Lemmas \ref{twisted1} and \ref{twisted2}, we see that if $\det \Phi (1-x_j)$ and $\det \Phi (1-x_j)$ are non-zero Laurent polynomials, then
$$
\frac{\det M^I_j}{\det \Phi (1-x_j)} = \pm \frac{\det M^I_k}{\det \Phi (1-x_k)} \in R(z),
$$
where $R(z)$ is the rational function field in $z$ over $R$.
If the dimension of the representation $\rho$ is even, then the sign in this formula is $+$.

Note that the Laurent polynomial ring $R[z^{\pm1}]$ over a unique factorization domain $R$ is again a unique factorization domain.
Therefore, we take the greatest common divisor $\gcd_I (\det M^I_j)$ of $\det M^I_j$ with respect to the choice of the indices $I$.
This Laurent polynomial is well-defined up to a factor $\varepsilon z^c \ (\varepsilon \in R^\times, c \in \mathbb{Z})$.

\begin{cor} [{\cite[Corollaly 5]{Wada}}]
If $\det \Phi (1-x_j)$ and $\det \Phi (1-x_j)$ are non-zero, then
$$
\frac{\gcd_I (\det M^I_j)}{\det \Phi (1-x_j)} = \varepsilon z^c \frac{\gcd_I (\det M^I_k)}{\det \Phi (1-x_k)}\ \ (\varepsilon \in R^\times, c \in \mathbb{Z}).
$$
\end{cor}

\begin{defi}
The \textbf{twisted Alexander invariant} $\Delta_{G, \rho}(z)$ of the group $G$ associated with the representation $\rho$ is defined as a rational expression
$$
\Delta_{G, \rho}(z):=\frac{\gcd_I (\det M^I_j)}{\det \Phi (1-x_j)}
$$
provided  $\det \Phi (1-x_j) \neq 0$.
If $m < l-1$, we define $\Delta_{G, \rho}(z) := 0$.
\end{defi}

Up to a factor $\varepsilon z^c \ (\varepsilon \in R^\times, c \in \mathbb{Z})$, the twisted Alexander invariant $\Delta_{G, \rho}(z)$ is well-defined independent of the choice of the number $j$.
Moreover, this is an invariant of the group $G$, the associated homomorphism $\alpha$, and the representation $\rho$.
In other words, the twisted Alexander invariant $\Delta_{G, \rho}(z)$ is independent of the choice of the presentation of $G$.
Further, if two representations $\rho_1$ and $\rho_2$ are equivalent, then $\Delta_{G, \rho_1}(z) \equiv \Delta_{G, \rho_2}(z)$.

\section{Proof of the main theorem}
\subsection{The case of $n=3$} \label{n=3}
In this case, $B_3$ is given by
$$
\left\langle \sigma_1, \sigma_2 \relmiddle| \sigma_1 \sigma_2 \sigma_1 = \sigma_2 \sigma_1 \sigma_2 \right\rangle.
$$
We write $r = \sigma_1 \sigma_2 \sigma_1 \sigma_2^{-1} \sigma_1^{-1} \sigma_2^{-1}$.
Then the Fox derivations of $r$ are given by
$$
\frac{\partial r}{\partial \sigma_1} = 1 + \sigma_1 \sigma_2 - \sigma_1 \sigma_2 \sigma_1 \sigma_2^{-1} \sigma_1^{-1}
$$
and
$$
\frac{\partial r}{\partial \sigma_2} = \sigma_1 - \sigma_1 \sigma_2 \sigma_1 \sigma_2^{-1} - \sigma_1 \sigma_2 \sigma_1 \sigma_2^{-1} \sigma_1^{-1} \sigma_2^{-1}.
$$
Thus
$$
\widetilde{\phi} \left(\frac{\partial r}{\partial \sigma_1}\right) = 1 + \sigma_1 \sigma_2 - \sigma_2
\ \ \ \textrm{and}\ \ \ 
\widetilde{\phi} \left(\frac{\partial r}{\partial \sigma_2}\right) = \sigma_1 - \sigma_2 \sigma_1 - 1.
$$
As the surjective homomorphism, we take the abelianization $\alpha : B_3 \longrightarrow \langle z \rangle$ given by
$$
\alpha(\sigma_1) = \alpha(\sigma_2) = z.
$$
Moreover, the Tong-Yang-Ma representation $TYM_3 : B_3 \longrightarrow GL_{3}(\mathbb{Z}[t^{\pm1}])$ is given by
$$
TYM_3(\sigma_1)=
\left(
\begin{array}{ccc}
0 & 1 & 0 \\
t & 0 & 0 \\
0 & 0 & 1
\end{array}
\right)
,\ \ 
TYM_3(\sigma_2)=
\left(
\begin{array}{ccc}
1 & 0 & 0 \\
0 & 0 & 1 \\
0 & t & 0
\end{array}
\right).
$$
The corresponding Alexander matrix $M$ is
$$
\left(
\Phi \left( \frac{\partial r}{\partial \sigma_1} \right) \ \Phi \left( \frac{\partial r}{\partial \sigma_2} \right)
\right)
=
\left(
\begin{array}{cccccc}
1-z & 0 & z^2 & -1 & z-z^2 & 0 \\ 
tz^2 & 1 & -z & tz & -1 & -z^2 \\
0 & -tz+tz^2 & 1 & -t^2z^2 & 0 & -1+z
\end{array}
\right)
$$
and hence
\begin{align*}
\det M_2 &=
\det
\left(
\begin{array}{ccc}
1-z & 0 & z^2 \\ 
tz^2 & 1 & -z \\
0 & -tz+tz^2 & 1
\end{array}
\right) \\
&= (1-z)(1-tz^2)(1+tz^3).
\end{align*}
We also have
\begin{align*}
\det \Phi (1 - \sigma_2) &=
\det
\left(
\begin{array}{ccc}
1-z & 0 & 0 \\
0 & 1 & -z \\
0 & -tz & 1
\end{array}
\right) \\
&= (1-z)(1-tz^2).
\end{align*}
Therefore, the twisted Alexander invariant of $B_3$ associated with $TYM_3$ is
\begin{align*}
\Delta_{B_3, TYM_3}(z)
&= \frac{\det M_2}{\det \Phi (1 - \sigma_2)} \\
&= 1+tz^3.
\end{align*}

\subsection{The case of $n \geq 4$}
The way of computations is based on \cite{Morifuji}.

We denote the relation between the generators $\sigma_i$ and $\sigma_j$ by $[i, j]$ and assign a number to each relation as follows:
$$
\begin{array}{llll}
r_1 = [1,2], & & &  \\
r_2 = [1,3], & r_3=[2,3], & &  \\
r_4 = [1,4], & r_5=[2,4], & r_6 = [3,4], &  \\
r_7 = [1,5], & r_8=[2,5], & r_9 = [3,5], & r_{10} = [4,5]  \\
\hspace{10pt} \vdots & & &\hspace{10pt} \ddots
\end{array}
$$

First we compute a denominator of the twisted Alexander invariant $\Delta_{B_n, TYM_n}(z)$.
\begin{lemma} \label{deno}
$
\det \Phi (1 - \sigma_{n-1}) = (1-z)^{n-2}(1-tz^2).
$
\end{lemma}

\begin{proof}
By a direct computation,
$$
\det  \Phi (1 - \sigma_{n-1})=
\det
\left(
\begin{array}{ccccc}
1-z & & & &  \\
 & \ddots & & &  \\
 & & 1-z & &  \\
 & & & 1 & -z \\
 & & & -tz & 1
\end{array}
\right)
=(1-z)^{n-2}(1-tz^2).
$$
\end{proof}

Next, we observe a numerator of $\Delta_{B_n, TYM_n}(z)$.

\begin{lemma} \label{num}
For any choice of the indices $I$, $\det M^I_{n-1}$ is divided by $ (1-z)^{n-2}(1-tz^2)$.
\end{lemma}

\begin{proof}
For each generator $\sigma_j$ and relation $r_i$,
$$
\widetilde{\phi} \left(\frac{\partial r_i}{\partial \sigma_j}\right) =
\begin{cases}
\sigma_k - 1 & (r_i = [k,j] , \ 1 \leq k \leq j-2) \\
\sigma_{j-1} - \sigma_j \sigma_{j-1} - 1 & (r_i = [j-1,j]) \\
1 + \sigma_{j} \sigma_{j+1} - \sigma_{j+1} & (r_i = [j,j+1]) \\
1 - \sigma_k & (r_i = [j,k] , \ j+2 \leq k \leq n-1) \\
0 & (\textrm{otherwise})
\end{cases},
$$
therefore
$$
\Phi \left( \frac{\partial r_i}{\partial \sigma_j} \right) =
\begin{cases}
(z-1) I_{k-1} \oplus
\left(
\begin{array}{cc}
-1 & z  \\
tz & -1 
\end{array}
\right)
\oplus (z-1) I_{n-k-1} \\
(-1+z-z^2) I_{j-2} \oplus
\left(
\begin{array}{ccc}
-1 & z-z^2 & 0 \\ 
tz & -1 & -z^2 \\
-t^2z^2 & 0 & -1+z
\end{array}
\right)
\oplus (-1+z-z^2) I_{n-j-1}  \\
(1-z+z^2) I_{j-1} \oplus
\left(
\begin{array}{ccc}
1-z & 0 & z^2 \\ 
tz^2 & 1 & -z \\
0 & -tz+tz^2 & 1
\end{array}
\right)
\oplus (1-z+z^2) I_{n-j-2} \\
(1-z) I_{k-1} \oplus
\left(
\begin{array}{cc}
1 & -z  \\
-tz & 1 
\end{array}
\right)
\oplus (1-z) I_{n-k-1} \\
0
\end{cases}.
$$

We denote the matrix $M_{n-1}$ by the column vectors ${\bm m_j} \ (1 \leq j \leq n-2)$:
$$
M_{n-1} = \left( \Phi \left( \frac{\partial r_i}{\partial \sigma_j} \right) \right)_{j \neq n-1}
= \left( {\bm m_1} , \ldots , {\bm m_{n-2}} \right)
$$
where
$$
{\bm m_j} = {^{t} \left( \Phi \left( \frac{\partial r_1}{\partial \sigma_j} \right) , \ldots , \Phi \left( \frac{\partial r_l}{\partial \sigma_j} \right) \right)}.
$$
Set $l$ to be the number of the relations of $B_n$, that is, $l=\displaystyle{\frac{(n-1)(n-2)}{2}}$.
Also, we denote the $i$-th column in ${\bm m_j}$ by $[i]_j$. Moreover, if we add $f$ times the column $[i]_j$ to the column $[k]_m$, denote it by
$$
[k]_m\  \xleftarrow{\textrm{add}}\  f \times [i]_j,
$$
where $f \in \mathbb{Z}[t^{\pm 1},z^{\pm 1}]$.

If we perform the following operations in order from the top:
\begin{align*}
[n-2]_{n-2} & \xleftarrow{\textrm{add}}  -tz \times [n-1]_{n-2} -z \times [n-3]_{n-3} + [n-2]_{n-3} \\ \relax
[n-3]_{n-3} & \xleftarrow{\textrm{add}}  -tz \times [n-2]_{n-3} -z \times [n-4]_{n-4} + [n-3]_{n-4} \\ \relax
&\hspace{10pt} \vdots \\ \relax
[3]_{3} & \xleftarrow{\textrm{add}}  -tz \times [4]_{3} -z \times [2]_{2} + [3]_{2} \\ \relax
[2]_{2} & \xleftarrow{\textrm{add}}  -tz \times [3]_{2} -z \times [1]_{1} + [2]_{1} \\ \relax
[1]_{1} & \xleftarrow{\textrm{add}}  -tz \times [2]_{1},
\end{align*}
then each column $[j]_j$ contains a common divisor $1-z$.
Therefore we take a term $1-z$ from the column $[j]_j \ (1 \leq j \leq n-2)$ as a divisor.
Hence we have $(1-z)^{n-2}$ as a common divisor of the matrix $M_{n-1}$.
We write $[j]'_j$ for the column $[j]_j$ divided by $1-z$. We replace the column $[j]_j$ in $M_{n-1}$ with $[j]'_j$ and denote the resulting matrix by $M'_{n-1}$.

Next if we perform the following operation:
\begin{align*}
[n-2]'_{n-2}\  \xleftarrow{\textrm{add}}\  -t \times [n]_{n-2}
& -t \times [n]_{n-3} -tz \times [n-1]_{n-3} + z \times [n-3]'_{n-3} \\
& -t \times [n]_{n-4} -tz \times [n-1]_{n-4} + z^2 \times [n-4]'_{n-4} \\
& \hspace{20pt} \vdots \\
& -t \times [n]_{2} -tz \times [n-1]_{2} + z^{n-4} \times [2]'_{2} \\
& -t \times [n]_{1} -tz \times [n-1]_{1} + z^{n-3} \times [1]'_{1},
\end{align*}
then the column $[n-2]'_{n-2}$ contains a common divisor $1-tz^2$.
Therefore we take a term $1-tz^2$ from this column as a divisor.
We write $\overline{[n-2]}_{n-2}$ for the column $[n-2]'_{n-2}$ divided by $1-tz^2$. We replace the column $[n-2]'_{n-2}$ in $M'_{n-1}$ with $\overline{[n-2]}_{n-2}$ and denote the resulting matrix by $\overline{M}_{n-1}$.
Accordingly we conclude that $\det M^I_{n-1}$ is divided by $(1-z)^{n-2}(1-tz^2)$ for any choice of indices $I$.
\end{proof}

In order to show $\gcd_I (\det M^I_{n-1}) = (1-z)^{n-2}(1-tz^2)$, we need the following lemma.

\begin{lemma} \label{only}
There exist the indices $I,I',I''$ such that\\
\begin{tabular}{rll}
{\rm (i)} & $\det \overline{M}^{I}_{n-1} = (1-tz^2)^{n-3}(1+tz^3)^{n-2}(1-z+z^2)^{(n-3)(n-2)}$ & $(n \geq 4)$, \\
{\rm (ii)} & $\det \overline{M}^{I'}_{n-1} = (1-tz^2)^{n-3}(1-z)^{(n-3)(n-2)}$ & $(n \geq 5)$, \\
{\rm (iii)} & $\det \overline{M}^{I''}_{n-1} = (1+tz^3)(1-z+z^2)^{n-2}(1-z)^{(n-4)(n-1)+1}$ & $(n \geq 4)$.
\end{tabular}
\end{lemma}

\begin{proof}
(i) For $n \geq 4$, let us consider an index $I$ corresponding to the $(n-2)$ row-blocks $[1,2], [2,3], \ldots, [n-2,n-1]$ in the matrix $M_{n-1}$.
Since
$$
\det \Phi (1 - \sigma_{j+1} + \sigma_{j} \sigma_{j+1}) = 
(1-z)(1-tz^2)(1+tz^3)(1-z+z^2)^{n-3},
$$
we have
\begin{align*}
\det M^I_{n-1}
&= \prod_{j=1}^{n-2} \textrm{det} \Phi (1 - \sigma_{j+1} + \sigma_{j} \sigma_{j+1}) \\
&= \left \{ (1-z)(1-tz^2)(1+tz^3) \}^{n-2} \right (1-z+z^2)^{(n-3)(n-2)}.
\end{align*}
Hence
$$
\det \overline{M}^I_{n-1} = (1-tz^2)^{n-3}(1+tz^3)^{n-2}(1-z+z^2)^{(n-3)(n-2)}
$$
for $n \geq 4$.

(ii) For $n \geq 5$, if we choose an index $I'$ corresponding to the $(n-2)$ row-blocks $[1,n-1], [2,n-1], \ldots , [n-3,n-1], [1,n-2]$, the diagonal blocks of square matrix $M^{I'}_{n-1}$ consist of $(n-3)$ matrices $\Phi (1 - \sigma_{n-1})$ and one $\Phi (\sigma_1 - 1)$.
Therefore
\begin{align*}
\det M^{I'}_{n-1}
&= \left \{ (1-tz^2)(1-z)^{n-2} \}^{n-3} \right (1-tz^2)(-1+z)^{n-2} \\
&= \pm (1-tz^2)^{n-2}(1-z)^{(n-2)(n-2)},
\end{align*}
where we used Lemma \ref{deno}.
Hence
$$
\det \overline{M}^{I'}_{n-1} = (1-tz^2)^{n-3}(1-z)^{(n-3)(n-2)}
$$
for $n \geq 5$.

(iii) In the matrix $M_{n-1}$ for $n \geq 4$, we replace the last row in $[k,n-1]$ with the last row in $[k,n-2] \ (1 \leq k \leq n-3)$.
Then the block $\Phi (1 - \sigma_{n-1})$ of each row-block $[k,n-1]$ changes to
$$
\begin{array}{ll}
(1-z) I_{n-2} \oplus
\left(
\begin{array}{cc}
1 & -z  \\
0 & 1-z 
\end{array}
\right) & (1 \leq k \leq n-4) \ \ \textrm{and} \\
(1-z) I_{n-2} \oplus
\left(
\begin{array}{cc}
1 & -z  \\
0 & 1-z+z^2 
\end{array}
\right) & (k=n-3).
\end{array}
$$
If we choose an index $I''$ corresponding to $n-3$ resulting row-blocks $[1,n-1], [2,n-1], \ldots, [n-3,n-1]$ and $[n-2,n-1]$, we have
\begin{align*}
\det M^{I''}_{n-1}
&= (1-z)^{(n-1)(n-4)} (1-z)^{n-2}(1-z+z^2) (1-z)(1-tz^2)(1+tz^3)(1-z+z^2)^{n-3} \\
&= (1-tz^2)(1+tz^3)(1-z+z^2)^{n-2} (1-z)^{(n-1)(n-3)}.
\end{align*}
Hence
$$
\det \overline{M}^{I''}_{n-1} = (1+tz^3)(1-z+z^2)^{n-2}(1-z)^{(n-4)(n-1)+1}
$$
for $n \geq 4$.
\end{proof}

By Lemmas \ref{num} and \ref{only}, for $n \geq 5$, $\gcd_I (\det M^I_{n-1}) = (1-z)^{n-2}(1-tz^2)$ holds.
In the case $n=4$, the matrix $M_3$ is as follows:
$$
\left(
\begin{array}{cccccccc}
1-z & 0 & z^2 & 0 & -1 & z-z^2 & 0 & 0 \\
t z^2 & 1 & -z & 0 & t z & -1 & -z^2 & 0 \\
0 & t z^2-t z & 1 & 0 & -t^2 z^2 & 0 & z-1 & 0 \\
0 & 0 & 0 & z^2-z+1 & 0 & 0 & 0 & -z^2+z-1 \\
1-z & 0 & 0 & 0 & 0 & 0 & 0 & 0 \\
0 & 1-z & 0 & 0 & 0 & 0 & 0 & 0 \\
0 & 0 & 1 & -z & 0 & 0 & 0 & 0 \\
0 & 0 & -t z & 1 & 0 & 0 & 0 & 0 \\
0 & 0 & 0 & 0 & z^2-z+1 & 0 & 0 & 0 \\
0 & 0 & 0 & 0 & 0 & 1-z & 0 & z^2 \\
0 & 0 & 0 & 0 & 0 & t z^2 & 1 & -z \\
0 & 0 & 0 & 0 & 0 & 0 & t z^2-t z & 1 \\
\end{array}
\right).
$$
If we choose the index $I=(1,2,5,6,7,8,11,12)$, then
$$
\det M^I_3 = (1-z)^2(1-tz^2)^2(1-tz^2+2tz^3-tz^4).
$$
Hence
$$
\det \overline{M}^I_3 = (1-tz^2)(1-tz^2+2tz^3-tz^4),
$$
and by Lemma \ref{only} (i) and (iii), we obtain the same result.
Therefore, by Lemma \ref{deno} and the above, Theorem \ref{twisted_TYM} immediately follows.

\section{Welded version}
\subsection{Definitions}
The welded braid group $WB_n$ is an extension of the braid group.
The braid group $B_n$ is interpreted as the fundamental group of the configuration space of $n$ points in the plane $\mathbb{R}^2$.
On the other hand, $WB_n$ is interpreted as the fundamental group of the configuration space of $n$ Euclidean, unordered, disjoint, unlinked circles in the $3$-ball $D^3$ lying on planes parallel to a fixed one (see \cite{Damiani} for example).
The welded braid group $WB_n$ has a presentation with generators $\{ \sigma_i,\tau_i \}_{i = 1,\ldots ,n-1}$ together with relations:
$$
\begin{cases}
\sigma_i \sigma_j = \sigma_j \sigma_i & (|i-j| \geq 2) \\ 
\sigma_i \sigma_{i+1} \sigma_i = \sigma_{i+1} \sigma_i \sigma_{i+1} & (i = 1, \ldots ,n-2) \\
\tau_i \tau_j = \tau_j \tau_i & (|i-j| \geq 2) \\
\tau_i \tau_{i+1} \tau_i = \tau_{i+1} \tau_{i} \tau_{i+1} & (i = 1, \ldots ,n-2) \\
\tau_i^2 = 1 &  (i = 1,\cdots ,n-1) \\
\sigma_i \tau_j = \tau_j \sigma_i & (|i-j| \geq 2) \\ 
\sigma_i \tau_{i+1} \tau_i = \tau_{i+1} \tau_i \sigma_{i+1} & (i = 1, \ldots ,n-2) \\
\tau_i \sigma_{i+1} \sigma_{i} = \sigma_{i+1} \sigma_{i} \tau_{i+1} & (i = 1, \ldots ,n-2)
\end{cases}
$$

\begin{figure}[h]
\begin{center}
\includegraphics[height=60pt]{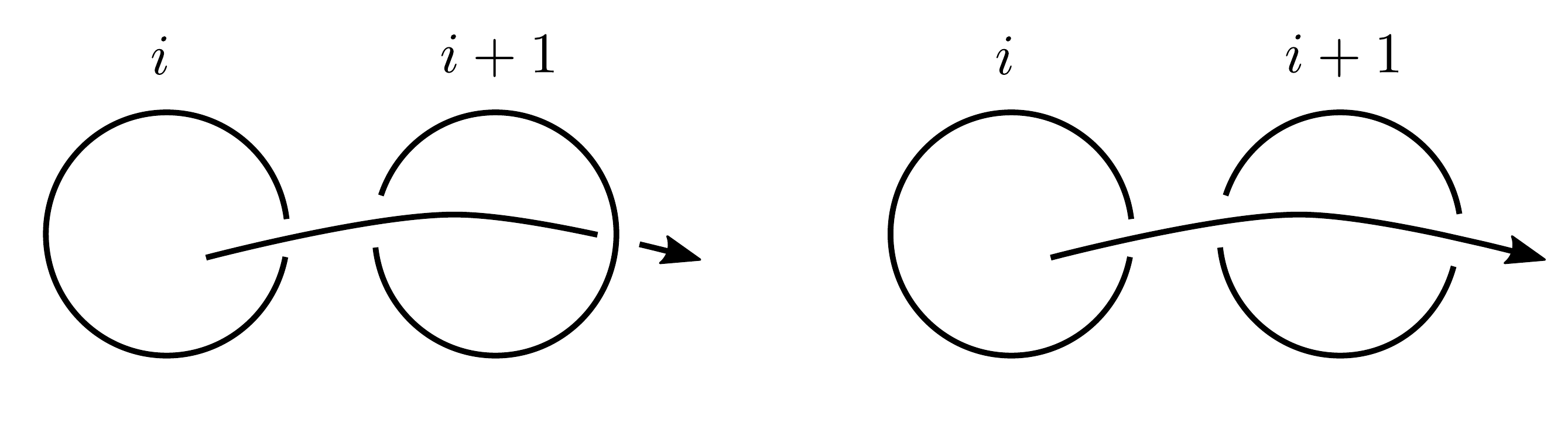}
\caption{Generators $\sigma_{i} , \tau_{i}$}
\label{welded_braid}
\end{center}
\end{figure}

\begin{figure}[h]
\begin{center}
\includegraphics[height=100pt]{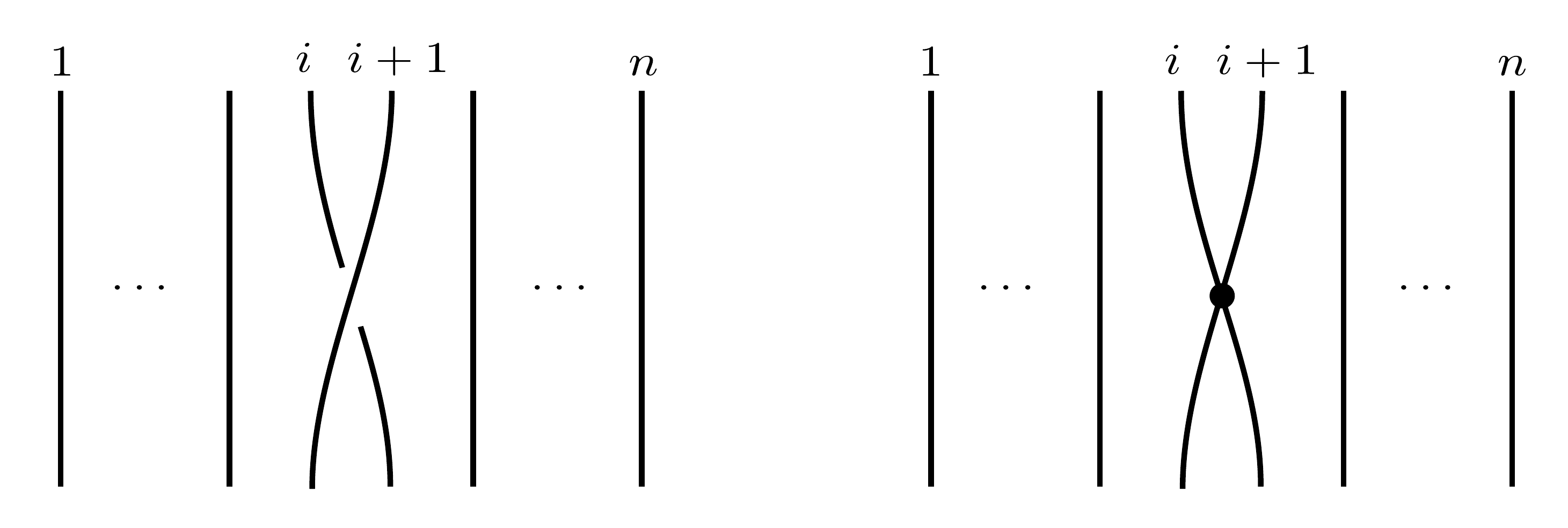}
\caption{Generators $\sigma_{i} , \tau_{i}$ as virtual braid diagrams}
\label{welded_diagram}
\end{center}
\end{figure}

The generators $\sigma_{i}$ is the loop permuting the $i$-th and the $(i+1)$-st circles by passing the $i$-th circle through the $(i+1)$-st circle, and $\tau_{i}$ permutes them by passing over (Figure \ref{welded_braid}).
Moreover, the element of $WB_n$ can be written as the virtual braid (Figure \ref{welded_diagram}). 

The Tong-Yang-Ma representation is exteded to the representation of the welded braid group in \cite{Bellingeri-Soulie}. This is given by
$$
\sigma_i \longmapsto
I_{i-1} \oplus
\left(
\begin{array}{cc}
0 & 1 \\
t & 0
\end{array}
\right)
\oplus I_{n-i-1}
\ \ {\rm and} \ \ 
\tau_i \longmapsto
I_{i-1} \oplus
\left(
\begin{array}{cc}
0 & a^{-1} \\
a & 0
\end{array}
\right)
\oplus I_{n-i-1}.
$$
We also call it the \textbf{Tong-Yang-Ma representation}, and denote it by w$TYM_n \colon WB_n \longrightarrow GL_{n} (\mathbb{Z} [t^{\pm1}, a^{\pm1}])$.
Now, we compute the twisted Alexander invariant of $WB_3$ associated with $\textrm{w}TYM_3$.

$WB_3$ is given by
$$
\left\langle \sigma_1, \sigma_2, \tau_1, \tau_2 \relmiddle|
\begin{array}{l}
\sigma_1 \sigma_2 \sigma_1 = \sigma_2 \sigma_1 \sigma_2, \tau_1 \tau_2 \tau_1 = \tau_2 \tau_1 \tau_1, \tau_1^2 = 1, \tau_2^2 = 1, \\
\sigma_1 \tau_2 \tau_1 = \tau_2 \tau_1 \sigma_2, \tau_1 \sigma_2 \sigma_1 = \sigma_2 \sigma_1 \tau_2
\end{array}
\right\rangle.
$$
We write
$$
\begin{cases}
r_1 = \sigma_1 \sigma_2 \sigma_1 \sigma_2^{-1} \sigma_1^{-1} \sigma_2^{-1}, \\
r_2 = \tau_1 \tau_2 \tau_1 \tau_2^{-1} \tau_1^{-1} \tau_1^{-1}, \\
r_3 = \tau_1^2, \\
r_4 = \tau_2^2, \\
r_5 = \sigma_1 \tau_2 \tau_1 \tau_2^{-1} \tau_1^{-1} \sigma_2^{-1}, \\
r_6 = \tau_1 \sigma_2 \sigma_1 \sigma_2^{-1} \sigma_1^{-1} \tau_2^{-1}.
\end{cases}
$$
Let $\alpha \colon WB_3 \longrightarrow \mathbb{Z} \cong \langle z \rangle$ be the surjective homomorphism given by 
$$
\alpha(\sigma_i) = z \ \ \ \textrm{and}\ \ \ \alpha(\tau_i) = 1
$$
for $i=1,2$.
Then the matrix $M_2$ obtained from the Alexander matrix $M$ by removing the second column, that is, the column corresponding to $\sigma_2$ is as follows:
$$
\left(
\begin{array}{ccccccccc}
 1-z & 0 & z^2 & 0 & 0 & 0 & 0 & 0 & 0 \\
 t z^2 & 1 & -z & 0 & 0 & 0 & 0 & 0 & 0 \\
 0 & t z^2-t z & 1 & 0 & 0 & 0 & 0 & 0 & 0 \\
 0 & 0 & 0 & 0 & 0 & a^{-2} & -1 & 0 & 0 \\
 0 & 0 & 0 & a & 1 & -a^{-1} & a & -1 & -a^{-1} \\
 0 & 0 & 0 & 0 & 0 & 1 & -a^2 & 0 & 0 \\
 0 & 0 & 0 & 1 & a^{-1} & 0 & 0 & 0 & 0 \\
 0 & 0 & 0 & a & 1 & 0 & 0 & 0 & 0 \\
 0 & 0 & 0 & 0 & 0 & 2 & 0 & 0 & 0 \\
 0 & 0 & 0 & 0 & 0 & 0 & 2 & 0 & 0 \\
 0 & 0 & 0 & 0 & 0 & 0 & 0 & 1 & a^{-1} \\
 0 & 0 & 0 & 0 & 0 & 0 & 0 & a & 1 \\
 1 & 0 & 0 & -1 & 0 & za^{-1} & -1 & z & 0 \\
 0 & 1 & 0 & t z & 0 & -a^{-1} & t z & -1 & 0 \\
 0 & 0 & 1 & 0 & a z-a & 0 & 0 & 0 & z-1 \\
 -z & 0 & za^{-1} & 1 & 0 & 0 & 0 & -z^2 & 0 \\
 a z & 0 & -z & 0 & 1 & 0 & 0 & 0 & -z^2 \\
 0 & 0 & 0 & 0 & 0 & 1 & -t^2 z^2 & 0 & 0 \\
\end{array}
\right).
$$

From Section \ref{n=3}, we obtain
$$
\det \Phi (1 - \sigma_2) = (1-z)(1-tz^2).
$$
Moreover by using Mathematica, the numerator of $\Delta_{WB_3,{\rm w}TYM_3}(z)$ is also $(1-z)(1-tz^2)$.
Therefore we obtain that $\Delta_{WB_3,{\rm w}TYM_3}(z) = 1$.
Then, from this fact and Theorem \ref{twisted_TYM},we naturally have the following question:

\begin{q}
For all $n \geq 3$, does $\Delta_{WB_n,{\rm w}TYM_n}(z)=1$ hold?
\end{q}

\section*{Acknowledgments}
The author would like to thank Takuya Sakasai for his careful reading of the paper and his helpful advice on this research.
Also, the author is grateful to the referees for careful reading and comments.

\bibliography{twisted}

\begin{thebibliography}{1}

\bibitem{Bellingeri-Soulie}
Paolo Bellingeri and Arthur Souli\'{e}.
\newblock A note on representations of welded braid groups.
\newblock {\em J. Knot Theory Ramifications}, 29(12):2050082, 21, 2020.

\bibitem{Damiani}
Celeste Damiani.
\newblock A journey through loop braid groups.
\newblock {\em Expo. Math.}, 35(3):252--285, 2017.

\bibitem{Morifuji}
Takayuki Morifuji.
\newblock Twisted {A}lexander polynomial for the braid group.
\newblock {\em Bull. Austral. Math. Soc.}, 64(1):1--13, 2001.

\bibitem{Tong-Yang-Ma}
Dian-Min Tong, Shan-De Yang, and Zhong-Qi Ma.
\newblock A new class of representations of braid groups.
\newblock {\em Comm. Theoret. Phys.}, 26(4):483--486, 1996.

\bibitem{Wada}
Masaaki Wada.
\newblock Twisted {A}lexander polynomial for finitely presentable groups.
\newblock {\em Topology}, 33(2):241--256, 1994.

\end{thebibliography}
\bibliographystyle{plain}

\end{document}